\documentclass[11pt, oneside]{article}   	
\usepackage{geometry}                		
\usepackage{graphicx}				
\usepackage{amssymb}

\newcommand{\N}{\mathcal{N}}

\title{Conjugacy of Cartan Subalgebras in Solvable Leibniz Algebras and Real Leibniz Algebras}
\author{Ernie Stitzinger and Ashley White}

\begin{document}
\maketitle

\noindent ABSTRACT. We extend conjugacy results from Lie algebras to their Leibniz algebra generalizations. The proofs in the Lie case depend on anti-commutativity. Thus it is necessary to find other paths in the Leibniz case. Some of these results involve Cartan subalgebras. Our  results can be used to extend other results on Cartan subalgebras. We show an example here and others will be shown in future work.\\

\noindent Keywords: Leibniz algebras, Conjugacy, Cartan subalgebras, maximal subalgebras

\noindent MSC 2010: 17A32\\

\noindent  A self-centralizing minimal ideal in a solvable Lie algebra is complemented and all complements are conjugate \cite{barnes2}. Such an algebra is called primitive. If the minimal ideal is not self-centralizing, then there is a bijection which assigns complements to their intersection with the centralizer \cite{stitz}. Such an intersection is the core of the complement. Towers extended this result to show that maximal subalgebras are conjugate exactly when their cores coincide, with the additional condition at characteristic $p$ that $L^2$ has nilpotency class less than $p$. Under this same condition, Cartan subalgebras are conjugate. We show these results hold for Leibniz algebras. The antisymmetry is used in the Lie algebra case, and so we amend the arguments in the Leibniz algebra case. We also show consequences for the conjugacy of Cartan subalgebras in solvable Lie algebras. We also consider conjugacy of Cartan subalgebras for real Leibniz algebras, obtaining an extension of a result of Barnes \cite{barnes2} in the Lie case.  Omirov \cite{omirov} treated conjugacy of  Cartan subalgebras in complex Leibniz algebras. For an introduction to Leibniz algebras, see \cite{demir}  \\

\noindent Let $A$ be an ideal in a solvable Leibniz algebra $L$. Also assume that at characteristic $p$, $A$ has nilpotency class $p$. As usual, for $a \in A$, let $\exp(L_a) = \sum_0^\infty \frac{1}{r!}(L_a)^r.$ Let $I(L, A)$ be the subgroup of the automorphism group of $L$ generated by $\exp L_a$. Two subalgeabras $U,V$ of $L$ are conjugate when there exists $\beta \in I(L, A)$ such that $\beta(U) =V$. If $A$ is abelian then $I(L, A)$ consists of all elements $I + L_a$. An algebra is primitive if it contains a self-centralizing minimal ideal, $A$. Barnes has shown in \cite{barnes1} that $A$ is complemented and all complements are conjugate under $I(L, A)$, which extends his Lie algebra result. The socle, Soc($L$), is the sum of all minimal ideals of $L$, and in the primitive case is equal to $A$. We record this as \\

\noindent {\bf Theorem 1} (\cite{barnes1}, 5.13) Let P be a primitive Leibniz algebra with socle $C$. Then $P$ splits over $C$ and all complements to $C$ in $P$ are conjugate under $I(L, C)$. \\  

\noindent We need the following result, which is the Leibniz algebra version of Lie algebra results found in \cite{barnes3}. \\

\noindent {\bf Lemma 1} Let $L$ be a solvable Leibniz algebra. If $\textrm{char}F = p$, also suppose that $L^2$ has nilpotency class less than $p$. Then,

\begin{itemize}
\item[$(i)$] for all $x \in L^2, \exp L_x$ exists,

\item[$(ii)$] if $U$ is a subalgebra of $L$, then every $\alpha \in I(U, U^2)$ has an extension $\alpha^* \in I(L, U^2)$,

\item[$(iii)$] if $A$ is an ideal of $L$, then every $\beta \in I(L/A, (L/A)^2)$ is induced by some $\beta^* \in I(L, L^2).$\\
\end{itemize}
\noindent{\bf Proof} $(i)$ For $x \in L^2$, $L_x(L) \subseteq L^2$. If $\textrm{char} F = 0$, then $L^2$ is nilpotent of class $c$ for some $c$, so $L_x^c(L^2) = 0$, and so $L_x^{c+1}(L) = 0.$  If $\textrm{char} F = p \neq 0$, then $L^2$ is nilpotent of class $c < p$, so $L_x^{p-1}(L^2) =0$, thus $L_x^p(L) = 0$.

$(ii)$ Let $\alpha \in I(U, U^2)$. Then $\alpha = \exp \lambda_{u_1} \exp \lambda_{u_2} \dots \exp \lambda_{u_n}$ for some $u_1 \dots u_n \in U^2$, where $\lambda_{u_i}$ is the restriction of $L_{u_i}$, the left multiplication operator, to $U$. Then $\alpha$ is just the restriction of $\alpha^* = \exp L_{u_1} \exp L_{u_2} \dots \exp L_{u_n} \in I(L, U^2)$.

$(iii)$ Let $\beta \in I(L/A, (L/A)^2)$. Then $\beta = \exp \lambda_{u_1} \exp \lambda_{u_2} \dots \exp \lambda_{u_n}$ , where $\lambda_{u_i}$ is left multiplication in (L/A) with $u_i \in (L/A)^2$. Since $(L/A)^2 = L^2 + A/A,$  $u_i = x_i + A$ for $x_i \in L^2$, and by $(i)$, $\exp L_{x_i}$ for $x_i \in L^2$ exists, so $\beta$ is induced by $\beta^* = \exp L_{x_1} \exp L_{x_2} \dots \exp L_{x_n} \in I(L, L^2)$.\\

\noindent We relax the condition that $A$ is self-centralizing to obtain the Leibniz algebra extension of the result in \cite{stitz}. We denote by $C_L(A)$ the centralizer of $A$ in $L$.\\

\noindent{\bf Theorem 2} Let $A$ be a minimal ideal of a solvable Leibniz algebra, $L.$ There exists a one to one correspondence between the distinct conjugate classes of complements to $A$ under $I(L, A)$, and the complements to $A$ in $C_L(A)$ that are ideals of $L$. \\

\noindent{\bf Proof} Let $M$ and $N$ be complements to $A$ in $L$ such that $M$ and $N$ are conjugate under $I(L, A).$ Define $M' = C_L(A) \cap M$ and $N' = C_L(A) \cap N$. Both $M'$ and $N'$ are complements to $A$ in $C_L(A)$, and are ideals of $L$. We show that $M' = N'$.  As $M$ and $N$ are conjugate, there exists a $\beta \in I(L,A)$ such that $N= \beta M$, and thus $N' = \beta M'$, but $M'$ is an ideal of $L$, so for $a \in A$, $M' = \exp L_a(M')$. Thus $M' = N'$.  Conversely let $M'$ be an ideal of $L$ and a complement to $A$ in $C_L(A)$. Then we have that $(A + M')/M'$ is its own centralizer in $L/M'$, thus is complemented in $L/M'$, and all complements are conjugate under $I(L/M', (A + M')/M')$. Let $M$ and $N$ be subalgebras of L that contain $M'$ such that $M/M'$ and $N/M'$ are complement to $(A+M')/M'$ in $L/M'.$ Then $M/M'$ and $N/M'$ are conjugate under $I(L/M', (A + M')/M')$, and by the correspondence theorem $M$ and $N$ are complements to $A$ in $L$. From Lemma 1, any $\beta^* \in I(L/M', (A+M')/M')$ is induced by some $\beta \in I(L, A)$, so $M$ and $N$ are conjugate in $I(L, A)$.\\

\noindent{\bf Corollary 1} Let $L$ be a solvable Leibniz algebra and $A$ be a minimal ideal of $L$ with complements $M$ and $N$. Then $M$ and $N$ are conjugate under $I(L, A)$ if and only if $M \cap C_L(A) = N \cap C_L(A)$.\\

\noindent The corollary follows immediately from Theorem 2. Towers extends these results further in \cite{towers} for the Lie case, we will extend them now for Leibniz algebras. For a subalgebra $U$ of $L$, its core, $U_L$ is the largest ideal of $L$ contained in $U$. $U$ is core-free if $U_L = 0$. Say two subalgebras are conjugate in $L$ if they are conjugate under $I(L, L) = I(L)$.\\

\noindent {\bf Lemma 2} Let $L$ be a solvable Leibniz algebra, and let $M, K$ be two core free maximal subalgebras of $L$. Then $M, K$ are conjugate under $\exp L_a = 1 + L_a$ for some $a \in L$, thus they are conjugate in $L$.\\

\noindent {\bf Proof} Let $A$ be a minimal abelian ideal of $L$. Then since $M, K$ are maximal subalgebras $L = A+M = A+K$. Note that $C_L(A) = A + (M \cap C_L(A))$, and $M \cap C_L(A)$ is an ideal of $L$. So $M \cap C_L(A) \subseteq M_L$, but $M$ is core-free, so $M \cap C_L(A) =0$. Thus $C_L(A) = A$, and the result follows from Theorem 1.  \\

\noindent{\bf Theorem 3} Let $L$ be a solvable Leibniz algebra over a field $F$. If $\textrm{char} F = p$, also suppose that $L^2$ has nilpotency class less than $p$. For maximal subalgebras $M,K$ of $L$, $M$ and $K$ are conjugate under $I(L, L^2)$ if and only if $M_L =  K_L$\\

\noindent {\bf Proof}  Suppose $M$ is conjugate to $K$ under $I(L, L^2)$. Then $M = \beta K$ where $\beta \in I(L, L^2)$, $\beta = \exp L_{x_1} \exp L_{x_2} \cdots \exp L_{x_n}$ for $x_i \in L^2$. Then $M_L = \beta K_L$, and further $\exp L_{x_i} K_L = K_L$ when $x_i \in L.$ So $\beta K_L = K_L$, thus $M_L = K_L$. Conversely let $M_L = K_L$. Then $M/M_L$ and $K/M_L$ are core-free maximal subalgebras of $L/M_L$, and are conjugate under $I(L/M_L, L/M_L)$ by Lemma 2, and thus also under $I(L/M_L, (L/M_L)^2)$. Then $M$ and $K$ are conjugate under $I(L, L^2)$ by Lemma 1. \\

\noindent A Cartan subalgebra of a Leibniz algebra $L$ is a nilpotent subalgebra $C$ such that $C = \N_L(C)$, where $\N_L(C)$ is the normalizer of $C$ in $L$. We now show that Cartan subalgebras of solvable Leibniz algebras over a field $F$ are conjugate in $L$, provided that when $\textrm{char}F = p$, $L^2$ has nilpotency class less than $p$. We will follow Barnes' Lie algebra approach from \cite{barnes1}, but make use of the following theorem.\\

\noindent{\bf Theorem 4} Let $L$ be a solvable Leibniz algebra, and let $M, K$ be complements of a minimal ideal $A$. If $M$ and $K$ are each Cartan subalgebras, then $M$ and $K$ are conjugate under $I(L, A)$.\\

\noindent{\bf Proof}  $M$ is a Cartan subalgebra and thus nilpotent, so consider the upper central series $Z_0(M) \subseteq Z_1(M) \subseteq \cdots \subseteq Z_j(M) = M$. Define $Z_j = Z_j(M) \cap C_L(M)$ to give the new series $0 \subseteq Z_0 \subseteq Z_1 \subseteq \cdots \subseteq Z_j = M \cap C_L(A)$. We have that $[M, Z_j] \subseteq Z_{j-1}, [A, Z_j]=0,$ and $[L, Z_j] \subseteq Z_{j-1}$. Suppose for some $j$, $Z_{j-1} \subseteq K$ and by way of contradiction $Z_j \nsubseteq K$. Then $[L, Z_j] \subseteq Z_{j-1} \subseteq K$, and $[K, Z_j] \subseteq Z_{j-1} \subseteq K$, so $Z_j \subseteq \N_L(K) = K$, contradiction. Thus $Z_j \subseteq K$ for all $j$, so $M \cap C_L(A) = K \cap C_L(A)$, and the result follows from Corollary 1.\\

\noindent{\bf Theorem 5}  Let $L$ be a solvable Leibniz algebra. If $\textrm{char}F = p \neq 0$, suppose further that $L^2$ is of nilpotency class less than $p$. Then the Cartan subalgebras of $L$ are conjugate under $I(L, L^2)$.\\

\noindent{\bf Proof} If $L$ is nilpotent the result holds. Suppose $L$ is not nilpotent. Let $H_1$ and $H_2$ be Cartan subalgebras of $L$ and let $A$ be a minimal ideal of $L$. Then $(H_1 + A)/A$ and $(H_2 + A)/A$ are conjugate under $I(L/A, (L/A)^2)$, thus $H_1 + A$ and $H_2 + A$ are conjugate under $I(L, L^2)$ by Lemma 1, so suppose that $H_1 + A = H_2 + A$. If $H_1 + A$ is a proper subalgebra of $L$, then $H_1$ and $H_2$ are conjugate under $I(H_1 + A, (H_1 + A)^2)$, and thus under $I(L, L^2)$ by Lemma 1.

Now suppose that $H_1 + A = H_2 + A = L$. $L$ is not nilpotent, so $H_1$ and $H_2$ are complements to $A$. By Theorem 4, $H_1$ and $H_2$ are conjugate under $I(L, A)$. Since $A$ is a minimal ideal in $L$, and thus an irreducible $L$-module, we have $al = -la$, or $al= 0$ for all $a \in A, l \in L$ \cite{barnes1}. Then $(LA)L = -L(LA)$ or $(LA)L = 0$, so $(LA)L \subseteq LA$, and $L(LA)  \subseteq LA$, so $LA$ is an ideal in $L$. But since $A$ is a minimal ideal, $LA = A$ or $LA = 0.$ $L$ is not nilpotent, so $LA \neq 0$. Thus $A = LA$, and as $LA$ is a subalgebra of $L^2$, $H_1$ and $H_2$ are conjugate under $I(L, L^2).$   \\

\noindent Now consider a Leibniz algebra $L$ over the field of real numbers, and the infinite series $\exp(Lx) = \sum_0^\infty \frac{1}{r!}(L_x)^r.$ For $x \in L$, the series converges and $\exp(L_x)$ is an automorphism of $L$. Call $I^*(L)$ the group of automorphisms of $L$ generated by $\exp(L_x)$ for $x \in L$. We find a criteria for two Cartan subalgebras of $L$ to be conjugate under $I^*(L)$.\\

\noindent{\bf Theorem 6} Let $L$ be a real Leibniz algebra with radical $R$. The conjugacy classes under $I^*(L)$ of Cartan subalgerbas of $L$ are in one to one correspondence with the conjugacy classes of $I^*(L/R)$ of Cartan subalgebras of $L/R$. Two Cartan subalgebras of $L$, $H_1, H_2$, are conjugate under $I^*(L) = I^*(L, L)$ if and only if $H_1 + R/R$ and $H_2 + R/R$ are conjugate under $I^*(L/R)$.\\

\noindent{\bf Proof} Assume $H_1$ and $H_2$ are conjugate under $I^*(L)$. Then $H_1 = \alpha H_2$ where $\alpha \in I^*(L).$ $\alpha$ induces an element, $\alpha^* \in I^*(L/R)$, so that $H_1+ R/R = \alpha^* (H_2 + R/R)$ and $H_1 + R/R$ and $H_2 + R/R$ are conjugate under $I^*(L/R).$ 

Conversely assume $H_1 + R/R$ and $H_2 + R/R$ are conjugate under $I^*(L/R)$ then $H_1 + R/R = \alpha^* (H_2 + R/R)$ for $\alpha^* \in I^*(L/R)$. $\alpha^*$ is induced by some $\alpha \in I^*(L)$, so $H_1+ R = \alpha ( H_2 + R)$, thus $H_1 + R$ and $H_2 + R$ are conjugate under $I^*(L)$. Then, say $H_1 + R = H_2 + R$. $H_1$ and $H_2$ are Cartan subalgebras of $H_1 + R$ and $H_2 + R$ respectively, so $H_1, H_2$ are Cartan subalgebras of $H_1 + R$ which is solvable, so $H_1$ and $H_2$ are conjugate under $I(H_1 + R)$ by Theorem 5. For any subalgebra $U$ of $L$, each element of $I^*(U)$ has an extension to an element in $I^*(L)$, so each element in $I(H_1 + R)$ has an extension to an element in $I^*(L)$, thus $H_1$ and $H_2$ are conjugate under $I^*(L).$    \\

\noindent The following conjugacy result is useful in extending results on Cartan subalgebras \cite{stitz1}  to their Leibniz algebra counterparts. \\

\noindent{\bf Theorem 7} Let L be a Leibniz algebra in which the intersection, J, of the terms in the lower central series is abelian. Let H and K be Cartan subalgebras of L. Then H and K are conjugate under  exp (L$_z$) where z $\in$ J. \\

\noindent{\bf Proof}  If J=0, then L is nilpotent and H=K. If J is a minimal ideal of L, then H and K are complements to J in L and the result holds by Theorem 4. If B is an ideal of L and is maximal as an ideal in J, then(H+B)/B and (K+B)/B are Cartan subalgebras of L/B and are conjugate under exp(L$_{(a+L/B)}$ in L/B by Theorem 4. Then H+B is conjugate to K+B in under exp(L$_a)$ in L where the image of H, H*, is a Cartan subalgebra of K+B. By induction, the Cartan subalgebras H* and K of K+B are conjugate under exp(L$_b)$, b$\in$ B. Then H and K are conjugate under exp(L$_{a+b})$ in L. Since a+b $\in$ J the result holds.\\

\noindent We now show the extension of two results from \cite{stitz1} to the Leibniz case, and will show further results in future work. Define the set $M(L)$ as the set of subalgebras, $A$, of $L$ such that: 
\begin{itemize}
\item[(1)]{$A$ can be joined to a $L$ by a chain of subalgebras of $L$ each maximal and self-normalizing in the next}
\item[(2)]{No maximal subalgebra of $A$ is self-normalizing in $A$}
\end{itemize} 

\noindent Notice that $M(L)$ is the set of subalgebras, $A$, of $L$ such that $A$ is nilpotent and $A$ satisfies condition $(1)$ above.  \\

\noindent Denote by $N(L)$ the nilradical of $L$. Say a Leibniz algebra is nilpotent of length $t$ if $t$ is the smallest natural number such that there exists a chain of ideals of $L$, $0 = L_0\subseteq L_1 \subseteq \cdots \subseteq L_t = L$ such that each factor is nilpotent. We show that when $L$ is nilpotent of length 2, the Cartan subalgebras of $L$ are precisely those subalgebras in $M(L)$, and when $L$ is nilpotent of length 3 each element of $M(L)$ is contained in exactly one Cartan subalgebra of $L$. We will need the following two lemmas, the proof of Lemma 4 is completely analogous to the Lie algebra case in Barnes \cite{barnes3}.  \\

\noindent {\bf Lemma 3} Let $L$ be a Leibniz algebra and $N$ an ideal of $L$, and $D$ be a subalgebra in $M(L)$, then $D+N/N \in M(L/N)$. \\

\noindent {\bf Proof} $D$ is nilpotent, so $D^m = 0$ for some $m$, and $(D + N)^m = D^m + N$, so $D^m + N /N =0$ and $D+N/N$ is nilpotent. There exists a chain of ideals $D \subseteq K_1 \subseteq K_2 \subseteq \cdots \subseteq K_n \subseteq L$, where $N_{K{i +1}}(K_i) = K_i$ and $K_i$ is maximal in $K_{i+1}$. Consider the new chain $D+N/N \subseteq K_1/N \subseteq K_2/N \subseteq \cdots \subseteq K_n/N \subseteq L/N$. Clearly $N_{K{i +1}/N}(K_i/N) = K_i/N$ and $K_i/N$ is maximal in $K_{i+1}/N$, so $D+N/N \in M(L/N).$\\

\noindent {\bf Lemma 4} Let $L$ be  solvable Leibniz algebra, $A$ an ideal in $L$ and $U$ a subalgebra such that $A \subseteq U \subseteq L$, where $U/A$ is a Cartan subalgebra of $L/A$. Then if $H$ is a Cartan subalgebra of $U$, $H$ is a Cartan subalgebra of $L$. \\

\noindent {\bf Theorem 8}  Let $L$ be a Leibniz algebra of nilpotent length 2 or less. Then the Cartan subalgebras of $L$ are those subalgebras in $M(L).$ \\

\noindent {\bf Proof} If $L$ is nilpotent the result holds, so assume $L$ is nilpotent of length 2 and proceed by induction on the dimension of $L$. Let $D \in M(L)$, then there exists a chain of subalgebras $0=D_0 \subseteq D_1 \subseteq D_2 \subseteq \cdots \subseteq D_t = L$, such that each $D_i$ is maximal and self-normalizing in the next. If $D$ is maximal in $L$, then $D = D_t$, thus $D$ is a Cartan subalgebra of $L$ and the result holds. Assume $D$ is not maximal in $L$, then by induction $D$ is a Cartan subalgebra of $D_t$, and as $D_t$ is maximal and self-normalizing in $L$, $L = D_t + N(L)$. Notice $L/N(L) = (D_t + N(L))/N(L)) \cong D_t/(D_t \cap N(L))$. $L/N(L)$ is nilpotent, and $D_t \subseteq \N_L(D_t \cap N(L))$, so $D_t/(D_t \cap N(L))$ is a Cartan subalgebra of $L/(D_t \cap N(L))$. $D$ is a Cartan subalgebra of $D_t$, so by Lemma 4, $D$ is a Cartan subalgebra of $L$.
Thus every element of $M(L)$ is a Cartan subalgebra of $L$.\\

\noindent {\bf Theorem 9} Let $L$ be nilpotent of length 3 or less, then each element of $M(L)$ is contained in exactly one Cartan subalgebra of $L$.\\

\noindent {\bf Proof} If $L$ is nilpotent of length less than three, then the result holds by the previous theorem, so let $L$ be nilpotent of length 3 and proceed by induction on the dimension of $L$. Let $D \in M(L)$, then $D + N(L)/N(L)$ is a Cartan subalgebra of $L/N(L).$ Let $E$ be a Cartan subalgebra of $D + N(L),$ then $D$ is a Cartan subalgebra of $L$ by Lemma 4. So $E + N(L) = D + N(L).$ Now let $S = E + N(L)$ and $S_x$ be the Fitting null component of $L_x$ acting on $S$. Define $K$ to be the intersection of all $S_x$ with $x \in D$. As $D$ in nilpotent, $D \subseteq K,$ we show $K$ is a Cartan subalgebra of $L$. Assume $y \in \N_L(K),$ but $y \notin K.$ Then $xy \in K,$ but $y \notin S_x,$ so $xy \notin S_x,$ contradiction. So $y \in K$ thus $K$ is its own normalizer in $L$. $K = K \cap (D + N(L)) = D + (K \cap N(L))$, and since $D$ is nilpotent by definition, and $(K \cap N(L))$ is a nilpotent ideal of $K$, $K$ is nilpotent by \cite{engel}.

Now assume $D$ is contained in two Cartan subalgebras of $L$, $E_1$ and $E_2.$ Then $E_1 + N(L) = D + N(L) = E_2 + N(L)$, so $E_1 = E_1 \cap (N(L) + D) = (E_1 \cap N(L)) + D$, and $E_2 = (E_2 \cap N(L)) + D.$ Suppose $A$ is a minimal ideal of $L$. Then $E_1 + A = E_2 +A$, and $E_1 + A$ is such that the intersection of its lower central series, $J$, is abelian, so by Theorem 7 all Cartan subalgebras of $E_1 + A$ are conjugate under $\exp(L_z) = (I + L_z), z \in J$. $E_1$ and $E_2$ are both Cartan subalgebras of $A + E_1$, so are conjugate under $I + L_z, z \in J$. So we have $E_1 \cap N(L) = (I + L_z)E_2 \cap N(L) = (I + L_z)(E_2 \cap N(L)) = E_2 \cap N(L)$, thus $E_1 = (E_1 \cap N(L)) + D = (E_2 \cap N(L)) + D = E_2$.

\end{document}